\documentclass{article}
\usepackage{amsmath,amsfonts,amstext}
\usepackage{geometry}                
\usepackage{graphicx} 
\usepackage{latexsym} 
\usepackage{pgf,pgfarrows,pgfnodes,pgfautomata,pgfheaps}
\usepackage{url}
\usepackage{amssymb}

\newtheorem{thm}{Theorem}[section]
\newtheorem{prop}[thm]{Proposition}

\newtheorem{lemma}[thm]{Lemma}

\newtheorem{cor}[thm]{Corollary}

\newtheorem{definitiontemp}[thm]{Definition}
\newenvironment{defn}{\begin{definitiontemp}
\normalfont}{\end{definitiontemp}}

\def\bec{\begin{cor}}
\def\enc{\end{cor}}

\def\bet{\begin{thm}}
\def\ent{\end{thm}}

\def\becor{\begin{cor}}
\def\encor{\end{cor}}

\def\bel{\begin{lem}}
\def\enl{\end{lem}}

\def\bedef{\begin{defn}}
\def\endef{\end{defn}}

\def\bep{\begin{prop}}
\def\enp{\end{prop}}

\newenvironment{pf}{\begin{trivlist}\item[\hskip\labelsep
{\it Proof.}]}{\end{trivlist}}

\newcommand{\set}[2]{\ensuremath{ \{ #1 : #2 \} }}

\newcommand{\Z}{\mathbb{Z}}

\newcommand{\Q}{\mathbb{Q}}

\newcommand{\R}{\mathcal{R}}
\newcommand{\C}{\mathcal{C}}
\newcommand{\A}{\mathcal{A}}
\newcommand{\B}{\mathcal{B}}
\newcommand{\D}{\mathcal{D}}
\newcommand{\E}{\mathcal{E}}
\renewcommand{\O}{\mathcal{O}}
\renewcommand{\P}{\mathbb{P}}
\renewcommand{\S}{\mathcal{S}}
\newcommand{\U}{\mathcal{U}}
\newcommand{\V}{\mathcal{V}}

\newcommand{\Wbar}{\overline{W}}

\newcommand{\Bbar}{\overline{\mathcal{B}}}

\newcommand{\hvec}{\vec{h}}

\newcommand{\xvec}{\vec{x}}
\newcommand{\Xvec}{\vec{X}}
\newcommand{\yvec}{\vec{y}}
\newcommand{\Yvec}{\vec{Y}}

\newcommand{\Zvec}{\vec{Z}}

\newcommand{\GL}[1]{\textbf{GL}_{#1}}

\newcommand{\comment}[1]{}

\newcommand{\la}{\langle}
\newcommand{\ra}{\rangle}

\newcommand{\cl}[1]{\text{cl}(#1)}
\newcommand{\Int}[1]{\text{Int}(#1)}

\newcommand{\SubQ}{\textbf{Sub}(\Q)}
\newcommand{\HTP}{\text{HTP}}
\newcommand{\HTPQ}{\text{HTP}(\Q)}

\def\diverges{\!\uparrow}
\def\converges{\!\downarrow}
\newcommand{\at}{\char'100}

\newcommand{\qed}{\hbox to 0pt{}\nobreak\hfill\rule{2mm}{2mm}}

\def\bfz{\boldsymbol{0}}
\def\s01{\ensuremath{\Sigma^0_1}}
\def\d02{\ensuremath{\Delta^0_2}}
\def\phi{\varphi}

\def\res{\!\!\upharpoonright\!}

\begin{document}

\title{Baire category theory and Hilbert's Tenth Problem inside $\Q$}


\author{Russell Miller}

\maketitle

\begin{abstract}
For a ring $R$, Hilbert's Tenth Problem HTP$(R)$ is the set of polynomial
equations over $R$, in several variables, with solutions in $R$.
We consider computability of this set for subrings $R$ of the rationals.
Applying Baire category theory to these subrings,
which naturally form a topological space, relates their
sets HTP$(R)$ to the set HTP$(\mathbb{Q})$, whose decidability
remains an open question.  The main result is that, for an arbitrary
set $C$, HTP$(\mathbb{Q})$ computes $C$ if and only if the subrings $R$
for which HTP$(R)$ computes $C$ form a nonmeager class.
Similar results hold for $1$-reducibility, for admitting
a Diophantine model of $\Z$, and for existential definability of $\Z$.
\end{abstract}

\section{Introduction}
\label{sec:intro}

The original version of Hilbert's Tenth Problem demanded an algorithm
deciding which polynomial equations from $\Z[X_1,X_2,\ldots]$ have solutions
in integers.  In 1970, Matiyasevic \cite{M70} completed work by Davis, Putnam
and Robinson \cite{DPR61}, showing that no such algorithm exists.
In particular, these authors showed that there exists a $1$-reduction from the
Halting Problem $\emptyset'$ to the set of such equations with solutions,
by proving the existence of a single polynomial $h\in\Z[Y,\Xvec]$
such that, for each $n$ from the set $\omega$
of nonnegative integers, the polynomial $h(n,\Xvec)=0$ has a solution
in $\Z$ if and only if $n$ lies in $\emptyset'$.
Since the membership in the Halting Problem was known
to be undecidable, it followed that Hilbert's Tenth Problem was also undecidable.

One naturally generalizes this problem to all rings $R$, defining
Hilbert's Tenth Problem for $R$ to be the set
$$ \HTP(R) = \set{f\in R[\Xvec]}{(\exists r_1,\ldots,r_n\in R^{<\omega})~f(r_1,\ldots,r_n)=0}.$$
Here we will examine this problem for one particular class:  the subrings $R$
of the field $\Q$ of rational numbers.  Notice that in this situation,
deciding membership in $\HTP(R)$ reduces to the question of deciding
this membership just for polynomials from $\Z[\Xvec]$, since one
readily eliminates denominators from the coefficients of a polynomial.
So, for us, $\HTP(R)$ will always be a subset of $\Z[X_1,X_2,\ldots]$.

Subrings $R$ of $\Q$ correspond bijectively to subsets $W$ of the set $\P$
of all primes, via the map $W\mapsto\Z[\frac1p:p\in W]$.  We write $R_W$
for the subring $\Z[\frac1p:p\in W]$. 
In this article, we will move interchangeably between subsets of $\omega$
and subsets of $\P$, using the bijection mapping $n\in\omega$ to the $n$-th prime $p_n$,
starting with $p_0=2$.
For the most part, our sets will be subsets of $\P$, but Turing reductions and
jump operators and the like will all be applied to them in the standard way.
Likewise, sets of polynomials, such as $\HTP(R)$, will be viewed as subsets of $\omega$,
using a fixed computable bijection from $\omega$ onto $\Z[\Xvec] = \Z[X_0,X_1,\ldots]$.

We usually view subsets of $\P$ as paths through the tree $2^{<\P}$,
a complete binary tree whose nodes are the functions from initial segments
of the set $\P$ into the set $\{ 0,1\}$.  This allows us to introduce a topology on the space
$2^{\P}$ of paths through $2^{<\P}$, and thus on the space of all subrings of $\Q$.
Each basic open set $\U_\sigma$ in this topology is given by a node $\sigma$ on the tree:
$\U_\sigma = \set{W\subseteq\P}{\sigma\subset W}$,
where $\sigma\subset W$ denotes that when $W$ is viewed as a function
from $\P$ into the set $2=\{ 0,1\}$ (i.e., as an infinite binary sequence),
$\sigma$ is an initial segment of that sequence.
Also, we put a natural measure $\mu$ on the class $\SubQ$ of all subrings of $\Q$:
just transfer to $\SubQ$ the obvious Lebesgue measure on the power set $2^{\P}$ of $\P$.
Thus, if we imagine choosing a subring $R$ by flipping a fair coin (independently for
each prime $p$) to decide whether $\frac1p\in R$, the \emph{measure} of a subclass
$\mathcal S$ of $\SubQ$ is the probability that the resulting subring will lie in $\mathcal S$.
Here we will focus on Baire category theory rather than on measure theory, however,
as the former yields more useful results.  For questions and results regarding measure theory,
we refer the reader to Section \ref{sec:measuretheory} and to the forthcoming \cite{M16}.

For all $W\subseteq \P$, we have the Turing reductions
$$ W\oplus\HTPQ \leq_T \HTP(R_W) \leq_T W'.$$
Indeed, each of these two Turing reductions is a \emph{$1$-reduction}.
\comment{
\begin{defn}
\label{defn:1reducible}
A $1$-reduction from $A$ to $B$ (where $A,B\subseteq\omega$)
is a computable injective function $h$ with domain $\omega$ such that
$$ (\forall n\in\omega)~[n\in A\iff h(n)\in B].$$
If such a function exists, we say that $A$ is \emph{$1$-reducible to $B$},
written $A\leq_1 B$.
\end{defn}
}
For instance, the Turing reduction from $HTP(R_W)$ to $W'$
can be described by a computable injection which maps
each $f\in\Z[\Xvec]$ to the code number $h(f)$ of an oracle Turing program
which, on every input, searches for a solution $\xvec$ to $f=0$ in $\Q$
for which the primes dividing the denominators of the coordinates
in $\xvec$ all lie in the oracle set $W$.  The reduction from $\HTPQ$
to $\HTP(R_W)$ uses the fact that every element of $\Q$ is a quotient
of elements of $R_W$, so that $f(\Xvec)$ has a solution in $\Q$
if and only if $Y^d\cdot f(\frac{X_1}{Y},\ldots,\frac{X_n}{Y})$ has a solution
in $R_W$ with $Y > 0$.  The condition $Y>0$ is readily expressed
using the Four Squares Theorem.

\section{Useful Facts}
\label{sec:lemmas}

The topological space $2^{\P}$ of all paths through $2^{<\P}$,
which we treat as the space of all subrings of $\Q$,
is obviously homeomorphic to \emph{Cantor space},
the space $2^{\omega}$ of all paths through the complete
binary tree $2^{<\omega}$.  Hence this space satisfies
the property of Baire, that no nonempty open set is meager.
We recall the relevant definitions.  Here as before, $\overline{\A}$ represents
the complement of a subset $\A\subseteq 2^{\P}$, and we will write
$\cl{\A}$ for the topological closure of $\A$ and $\Int{\A}$ for its interior.

\begin{defn}
\label{defn:meager}
A subset $\B\subseteq 2^{\P}$ is said to be \emph{nowhere dense}
if its closure $\cl{\B}$ contains no nonempty open subset of $2^{\P}$.
In particular, every set $\U_\sigma$ with $\sigma\in 2^{<\P}$
must intersect $\Int{\overline{\B}}$, the interior of the complement of $\B$.

The union of countably many nowhere dense subsets of $2^{\omega}$
is called a \emph{meager} set, or a \emph{set of first category}.
Its complement is said to be \emph{comeager}.
\end{defn}

The term ``Baire space'' is often used to name the particular space
$\omega^{\omega}$.  This terminology is confusing, since ``Baire space''
is also sometimes used to describe any space satisfying the property
of Baire.  The compact space $2^{\P}$ is certainly not homeomorphic to
the noncompact space $\omega^{\omega}$, but the property of Baire holds in both:
nonempty sets cannot be meager.

All sets $W\subseteq\omega$ satisfy $W\oplus\emptyset'\leq_T W'$, and for certain $W$,
Turing-equivalence holds here.  Indeed, it is known that the class
$$ \GL1 = \set{W\in 2^{\omega}}{W'\equiv_T W\oplus \emptyset'}$$
is comeager, although its complement is nonempty.
In computability theory, elements of $\GL1$ are called
\emph{generalized-low$_1$ sets}.  The low sets -- i.e., those $W$
with $W'\leq_T\emptyset'$ -- clearly lie in $\GL1$.
\begin{lemma}[Folklore]
\label{lemma:uniform}
There exists a Turing functional $\Psi$ such that
$\set{W\subseteq\omega}{\Psi^{W\oplus\emptyset'}=\chi_{W'}}$ is comeager.
It follows that $\GL1$ is comeager.
\end{lemma}
\begin{pf}  
Consider the following oracle program $\Psi$
for computing $W'$ from $W\oplus\emptyset'$.  With this oracle, on input $e$,
the program searches for a string $\sigma\subseteq W$ such that either
\begin{enumerate}
\item
$(\exists s)~\Phi_{e,s}^{\sigma}(e)\converges$; or
\item
$(\forall \tau\supseteq\sigma)(\forall s)~\Phi_{e,s}^{\tau}(e)\diverges$.
\end{enumerate}
The program uses its $\emptyset'$ oracle to check the truth of these
two statements for each $\sigma\subseteq W$.  If it ever finds that (1) holds,
it concludes that $e\in W'$; while if it ever finds that (2) holds,
it concludes that $e\notin W'$.  Thus, $\Psi^{W\oplus\emptyset'}$
can only fail to compute $W'$ if there exists some $e\notin W'$
such that, for every $n$, some $\tau\supset W\res n$ has $\Phi_e^{\tau}(e)\converges$.
This can happen, but for each single $e$, the set of those $W$
for which this happens constitutes the boundary of the open set
$\set{W}{e\in W'}$.  This boundary is nowhere dense (cf.\ Lemma
\ref{lemma:boundarymeager} below), so the union of these sets
(over all $e$) is meager, and $\Psi^{W\oplus\emptyset'}=\chi_{W'}$
for every $W$ outside this meager set.
\qed\end{pf}
$\GL1$ also has measure $1$, but no single Turing functional computes $W'$
from $W\oplus\emptyset'$ uniformly on a set of measure $1$.

\begin{lemma}[Folklore]
\label{lemma:majorityvote}
If $A \not\geq_T B$, then the class $\C=\set{W}{A\oplus W\geq_T B}$
is meager.
\end{lemma}
\begin{pf}
To show that $\C$ is meager, define $\C_e=\set{W\subseteq\P}{\Phi_e^{A\oplus W}=\chi_B}$,
so $\C=\cup_e\C_e$.
We claim that, if $\sigma\in 2^{\P}$ and $\U_\sigma\subseteq\cl{\C_e}$,
the following hold.
\begin{enumerate}
\item
$\forall x\forall \tau\supseteq\sigma~[\Phi_e^{A\oplus\tau}(x)\diverges\text{~or~}
\Phi_e^{A\oplus\tau}(x)\converges=\chi_B(x)]$.
\item
$\forall x\exists\tau\supseteq\sigma~[\Phi_e^{A\oplus\tau}(x)\converges]$.
\end{enumerate}

To see that (1) holds, suppose $\Phi_e^{A\oplus\tau}(x)\converges$.
With $\U_{\tau}\subseteq \U_{\sigma}\subseteq \text{cl}(\C_e)$, some $W\in\C_e$
must have $\tau\subseteq W$.  But then $\chi_B(x)=\Phi_e^{A\oplus W}(x)\converges
=\Phi_e^{A\oplus\tau}(x)$.

To see (2), fix any $W\in\C_e$ with $\sigma\subseteq W$:  such a $W$ must exist,
since $\U_{\sigma}\subseteq\text{cl}(\C_e)$.  Then we can take $\tau$ to be the
restriction of this $W$ to the use of the computation $\Phi_e^{A\oplus W}(x)$
(or $\tau=\sigma$ if the use is $<|\sigma|$).

But now every $\C_e$ must be nowhere dense, since any $\sigma$ satisfying
(1) and (2) would let us compute $B$ from $A$:
given $x$, just search for some $\tau\supseteq\sigma$ and some $s$ for which
$\Phi_{e,s}^{A\oplus\tau}(x)\converges$.  By (2), our search would discover such a $\tau$
eventually, and by (1) we would know $\chi_B(x)=\Phi_{e,s}^{A\oplus\tau}(x)$.
Since $A \not\geq_T B$, this is impossible.
\qed\end{pf}

Finally, on a separate topic, it will be important for us to know that whenever
$R$ is a semilocal subring of $\Q$, we have $\HTP(R)\leq_1\HTPQ$.
Indeed, both the Turing reduction and the $1$-reduction are uniform in the complement.
(The result essentially follows from work of Julia Robinson in \cite{R49}.  For a proof by
Eisentr\"ager, Park, Shlapentokh, and the author, see \cite{EMPS15}.)
Recall that the \emph{semilocal} subrings of $\Q$ are precisely those of the form
$R_W$ where the set $W$ is cofinite in $\P$, containing all but finitely many primes.
\begin{prop}[see Proposition 5.4 in \cite{EMPS15}]
\label{prop:semilocal}
There exists a computable function $G$ such that
for every $n$, every finite set $A_0=\{ p_1,\ldots,p_n\}\subset\P$ and every $f\in\Z[\Xvec]$,
$$ f\in\HTP(R_{\P-A_0}) \iff G(f,\la p_1,\ldots,p_n\ra) \in \HTPQ.$$
That is, $\HTP(R_{\P-A_0})$ is $1$-reducible to $\HTPQ$ for all semilocal
$R_{\P-A_0}$, uniformly in $A_0$.
\qed\end{prop}
The proof in \cite{EMPS15}, using work from \cite{K16},
actually shows how to compute, for every prime $p$,
a polynomial $f_{p}(Z,X_1,X_2,X_3)$ such that for all rationals $q$, we have
$$ q\in R_{\P-\{p\}} \iff f_{p}(q,\Xvec)\in\HTPQ.$$
Therefore, an arbitrary $g(Z_0,\ldots,Z_n)$ has a solution in $R_{\P-A_0}$ if and only if
$$(g(\Zvec))^2 + \sum_{p\in A_0, j\leq n} (f_{p}(Z_j,X_{1j},X_{2j},X_{3j}))^2$$
has a solution in $\Q$.

\section{Baire Category and Turing Reducibility}
\label{sec:leq_T}

For a polynomial $f\in\Z[\Xvec]$ and a subring $R_W\subseteq\Q$,
there are three possibilities.  First, $f$ may lie in $\HTP(R_W)$.
If this holds for $R_W$, the reason is finitary:  $W$ contains a certain
finite (possibly empty) subset of primes generating the denominators
of a solution.  Second, there may be a finitary reason why $f\notin\HTP(R_W)$:
there may exist a finite subset $A_0$ of the complement $\Wbar$
such that $f$ has no solution in $R_{\P-A_0}$.  For each finite $A_0\subset\P$,
the set $\HTP(R_{\P-A_0})$ is $1$-reducible to $\HTPQ$, by Proposition \ref{prop:semilocal};
indeed the two sets are computably isomorphic, with a computable permutation
of $\Z[\Xvec]$ mapping one onto the other.  Therefore, the existence
of such a set $A_0$ (still for one fixed $f$) is a $\Sigma_1^{\HTPQ}$ problem.

The third possibility is that neither of the first two holds.  An example
is given in \cite{M16}, where it is shown that a particular polynomial $f$
fails to lie in $\HTP(R_{W_3})$, where $W_3$ is the set of all primes
congruent to $3$ modulo $4$, yet that, for every finite set $V_0$ of primes,
there exists some $W$ disjoint from $V_0$ with $f\in\HTP(R_W)$.
We consider sets such as this $W_3$ to be on the \emph{boundary} of $f$,
in consideration of the topology of the situation.
The set $\A(f)=\set{W}{f\in\HTP(R_W)}$ is open in the usual
topology on $2^\P$, since, for any solution of $f$ in $R_W$
and any $\sigma\subseteq W$ long enough to include
all primes dividing the denominators in that solution, every other
$V\supseteq\sigma$ will also contain that solution.  Moreover,
one can computably enumerate the collection of those $\sigma$ such that
the basic open set $\U_\sigma=\set{W}{\sigma\subseteq W}$
is contained within $\A(f)$.  The set $\Int{\overline{\A(f)}}$
is similarly a union of basic open sets, and these can be enumerated
by an $\HTPQ$-oracle, since $\HTPQ$ decides $\HTP(R)$
uniformly for every semilocal ring $R$.
The \emph{boundary} $\B(f)$ of $f$ remains:  it contains those $W$
which lie neither in $\A(f)$ nor in $\Int{\overline{\A(f)}}$.
The boundary can be empty, but need not be, as seen
in the example mentioned above.

It follows quickly from Baire category theory that
the boundary set for a polynomial $f\in\Z[\Xvec]$
must be nowhere dense.  In general the boundary set $\partial\A$ of a set $\A$
within a space $\S$ is defined to equal $(\S-\Int{\A}-\Int{\overline{\A}})$,
and thus is always closed.

\begin{lemma}
\label{lemma:boundarymeager}
For every open set $\A$ in a Baire space $\S$, the boundary set
$\partial\A$ is nowhere dense.  In particular,
for each $f\in\Z[\Xvec]$, 
the boundary set $\B(f) = \partial(\A(f))$
must be nowhere dense.  Hence the \emph{entire boundary set}
$$ \B = \set{W\subseteq\P}{(\exists f\in\Z[\Xvec])~W\in \B(f)} = \cup_{f\in\Z[\Xvec]} \B(f)$$
is meager.
\end{lemma}
\begin{pf}
Since $\A$ is open, every open subset $\V$ of the closure of $\partial\A$ (namely $\partial\A$ itself)
lies within the complement $\overline{\A}$, hence within $\Int{\overline{\A}}$,
which is also disjoint from $\partial\A$.  This proves that
$\partial\A$ is nowhere dense.
Hence $\B$, the countable union of such sets, is meager.
\qed\end{pf}

For a set $W$ to fail to lie in $\B$, it must be the case that for every polynomial
$f$, either $f\in\HTP(R_W)$ or else some finite initial segment of $W$
rules out all solutions to $f$.  This is an example of the concept of
\emph{genericity}, common in both computability and set theory, so we adopt the term here.
With this notion, we can show not only that $\HTP(R_W)\leq W\oplus\HTPQ$
for all $W$ in the comeager set $\Bbar$, but indeed that the reduction is uniform
on $\Bbar$.
\begin{defn}
\label{defn:HTP-generic}
A set $W\subseteq\P$ is \emph{HTP-generic} if $W\notin\B$.
In this case we will also call the corresponding subring $R_W$ HTP-generic.
By Lemma \ref{lemma:boundarymeager}, HTP-genericity is comeager.
\end{defn}

\begin{prop}
\label{prop:uniformcomeager}
There is a single Turing reduction $\Phi$ such that the set
$$ \set{W\subseteq\P}{\Phi^{W\oplus\HTPQ}=\chi_{\HTP(R_W)}}$$
is comeager.  Hence $\HTP(R_W)\equiv_T W\oplus\HTPQ$ for every HTP-generic set $W$.
\end{prop}
\begin{pf}
Given $f\in\Z[\Xvec]$ as input, the program for $\Phi$ simply searches
for either a solution $\xvec$ to $f=0$ in $\Q$ for which all primes
dividing the denominators lie in the oracle set $W$, or else a finite
set $A_0\subseteq \Wbar$ such that the $\HTPQ$ oracle, using Proposition
\ref{prop:semilocal}, confirms
that $f\notin\HTP(R_{\P-A_0})$.  When it finds either of these, it outputs
the corresponding answer about membership of $f$ in $\HTP(R_W)$.
If it never finds either, then $W\in \B(f)$, and so this process succeeds
for every $W$ except those in the meager set $\B$.
\qed\end{pf}

\begin{cor}
\label{cor:comeager}
For every set $C\subseteq\omega$, the following are equivalent
\begin{enumerate}
\item
$C\leq_T\HTPQ$.
\item
$\set{W\subseteq\P}{C\leq_T\HTP(R_W)}=2^{\P}$ .
\item
$\set{W\subseteq\P}{C\leq_T\HTP(R_W)}$ is comeager.
\item
$\set{W\subseteq\P}{C\leq_T\HTP(R_W)}$ is not meager.
\end{enumerate}
\end{cor}
This opens a new possible avenue to a proof of undecidability of $\HTPQ$:
one need not address $\Q$ itself, but only show that for most subrings $\R_W$,
$\HTP(R_W)$ can decide the halting problem (or some other fixed undecidable set $C$).
Constructions in the style of \cite[Theorem 1.3]{P03} offer an approach
to the problem along these lines:  that theorem, proven by Poonen, shows that the set of such
subrings has size continuum and is large in certain other senses,
although the set of subrings given there is nowhere dense and
therefore does not by itself enable us to apply Corollary \ref{cor:comeager}.
\begin{pf}
Trivially $(1\!\!\implies\! 2\!\!\implies\! 3)$, since all $W$ satisfy $\HTPQ\leq_T\HTP(R_W)$,
and $(3\implies 4)$ holds in every Baire space.
So assume (4).  Then by Proposition \ref{prop:uniformcomeager},
$C\leq_T W\oplus\HTPQ$ holds on a non-meager set, as the intersection
of a non-meager set with a comeager set cannot be meager.
So by Lemma \ref{lemma:majorityvote}, $C\leq_T\HTPQ$.
\qed\end{pf}



Of course, $\HTP(R_W)$ always computes $W$, and so
$\HTP(R_W)$ is undecidable whenever $W$ is not computable.
However, we can strengthen the above statements a little further.
\begin{prop}
\label{prop:Delta2}
If $\HTPQ$ is decidable, then the following class is meager:
$$ \D=\set{W\subseteq\P}{(\exists D\leq_T\emptyset')~\emptyset <_T D\leq_T\HTP(R_W)}.$$
\end{prop}
Thus, while undecidability of $\HTP(R_W)$ is a given whenever $W>_T\emptyset$, the ability of
$\HTP(R_W)$ to compute any noncomputable $\Delta^0_2$ set $D$ is of real interest.
Even if different subrings in this class compute many distinct such sets $D$
-- and even if these sets all form \emph{minimal pairs}, i.e., their degrees
all have pairwise infimum $\bfz$ --
there are only countably many such $D$, which is the key to the proof.
\begin{pf}
Let $\la D_n\ra_{n\in\omega}$ be any (necessarily noneffective)
enumeration of the noncomputable $\Delta^0_2$ sets, and define
$$\D_n = \set{W\subseteq\P}{D_n\leq_T\HTP(R_W)}.$$
If a single $\D_n$ were not meager, then the intersection $(\D_n\cap\overline{\B})$
would also not be meager, since the entire boundary set $\B$ is meager.
But every $W\in\Bbar$ has $\HTP(R_W)\leq_T W\oplus\HTPQ$,
and so $\set{W\subseteq\P}{D_n\leq W\oplus\HTPQ}$ would also
not be meager.  By Lemma \ref{lemma:majorityvote},
this would imply $D_n\leq_T\HTPQ$.  Therefore, the assumption
that $\HTPQ$ is decidable ensures that every $\D_n$ is meager,
making their countable union $\D$ is meager as well.
\qed\end{pf}

\section{$\bf{1}$-Reducibility and Baire Category}
\label{sec:1Baire}

In Section \ref{sec:leq_T} we examined classes of subsets of $\P$
defined by Turing reductions involving $\HTP(R_W)$.
Here we replace Turing reducibility by $1$-reducibility
and ask similar questions about classes so defined.
It is not known whether there exists a subring $R\subseteq\Q$
for which $\emptyset'\leq_T\HTP(R_W)$ but $\emptyset'\not\leq_1\HTP(R_W)$,
and we have no good candidates for such a subring.
Ever since the original proof of undecidability of Hilbert's
Tenth Problem in \cite{DPR61,M70}, every Turing reduction
ever given from the Halting Problem to any $\HTP(R)$
with $R\subseteq\Q$ has in fact been a $1$-reduction.
Of course, if $\emptyset'\leq_1\HTPQ$, then $\emptyset'\leq_1\HTP(R)$
for all subrings $R$, so in some sense $\Q$ itself is the ``only'' candidate.

We have a result for $1$-reducibility analogous to Corollary \ref{cor:comeager},
but the proof is somewhat different.

\begin{thm}
\label{thm:1Baire}
For every set $C\subseteq\omega$ with $C\not\leq_1\HTPQ$, the following
class is meager:
$$ \O = \set{W\subseteq\P}{C\leq_1\HTP(R_W)}.$$
\end{thm}
\begin{pf}
One naturally views $\O$ as the union of countably many
subclasses $\O_e$, where
$$ \O_e = \set{W\subseteq\P}{C\leq_1\HTP(R_W)
\text{~via~}\phi_e}.$$
Of course, for those $e$ for which the $e$-th Turing function
$\phi_e$ is not total, this class is empty.  We claim that
if any one of these $\O_e$ fails to be nowhere dense,
then $C\leq_1\HTPQ$, contrary to the assumption
of the theorem.

Suppose that indeed $\O_e$ fails to be nowhere dense,
and fix a $\sigma$ for which $\U_\sigma\subseteq\cl{\O_e}$.
Let $A_0=\sigma^{-1}(0)$ contain those primes excluded
from all $W\in\U_\sigma$, and set $R=R_{(\P-A_0)}$.
Now whenever $n\in C$ and $W\in\O_e$,
the polynomial $\phi_e(n)$ must lie in $\HTP(R_W)$.
Since some $W\in\O_e$ lies in $U_\sigma$, we must have
$\phi_e(n)\in\HTP(R)$, because $R_W\subseteq R$
whenever $W\in\U_\sigma$.  On the other hand, suppose $n\notin C$.
If $R$ contained a solution to the polynomial $\phi_e(n)$,
then some $\tau\supseteq\sigma$ would by itself invert
the finitely many primes required to generate this solution,
and thus we would have $\U_\tau\cap\O_e=\emptyset$.
With $\U_\sigma\subseteq\cl{\O_e}$, this is impossible,
and so, whenever $n\notin C$, we have
$\phi_e(n)\notin\HTP(R)$.

Thus $R$ itself lies in $\O_e$, as $\phi_e$ is a $1$-reduction
from $C$ to $\HTP(R)$.  But $R$ is semilocal,
inverting all primes $p$ except those with $\sigma(p)=0$.
By Proposition \ref{prop:semilocal}, we have
$\HTP(R)\leq_1 \HTPQ$, and so $C\leq_1\HTPQ$.
\qed\end{pf}

\section{More All-Or-Nothing Laws}
\label{sec:diophantine}

This section proves two similar results, one about subrings of $\Q$
which admit diophantine models and one about subrings which admit
existential definitions of the integers within the subring.  In both cases,
the result is a sort of zero-one law:  that the given phenomenon must
either hold almost everywhere (i.e., on a comeager set of subrings)
or almost nowhere (i.e., on a meager set).  We begin with the diophantine models.

\begin{defn}
\label{defn:diophantine}
In a ring $R$, a \emph{diophantine model of $\Z$}
consists of three polynomials $h$, $h_+$, and $h_{\times}$,
with $h\in R[X_1,\ldots,X_n,\Yvec]$
and $h_+,h_{\times}\in R[X_1,\ldots,X_{3n},\Yvec]$ (for some $n$),
such that the set
$$\set{\xvec\in R^n}{(\exists\yvec\in R^{<\omega})~h(\xvec,\yvec)=0}$$
(equivalently, $\set{\xvec\in R^n}{h(\xvec,\Yvec)\in\HTP(R)}$)
is isomorphic to the structure $(\Z,+,\cdot)$
under the binary operations whose graphs are defined by
$$ \set{(\xvec_1,\xvec_2,\xvec_3)\in R^{3n}}{h_+(\xvec_1,\xvec_2,\xvec_3,\Yvec)\in\HTP(R)}$$
for addition and the corresponding set with $h_{\times}$ for multiplication.
\end{defn}

If a computable ring $R$ admits a diophantine model of $\Z$, then
$\HTP(\Z)$ can be coded into $\HTP(R)$, and so $\emptyset'\equiv_1\HTP(\Z)\leq_1\HTP(R)$.
For subrings $R_W$ of $\Q$ for which $\emptyset'\not\leq_T W$,
this is the only known method of showing that $\emptyset'\leq_T\HTP(R_W)$
(apart from the original proof by Matiyasevich, Davis, Putnam, and Robinson
for the case $W=\emptyset$, of course, which is what allows this method to succeed).

\begin{defn}
\label{defn:Dclass}
$\D$ is the class of subrings of $\Q$ admitting a diophantine model:
$$ \D = \set{W\subseteq\P}{R_W\text{~admits a diophantine model of~}\Z}.$$
\end{defn}

In this section we address the question of the size of the class $\D$.  The main
result fails to resolve this question, but shows it to have an ``all-or-nothing''
character.

\begin{thm}
\label{thm:diophantine}
One (and only one) of the following two possibilities holds.
\begin{enumerate}
\item
The class $\D$ is meager.
\item
There exists a particular triple $(h,h_+,h_{\times})$ of polynomials
over $\Z$ and a finite binary string $\sigma\in 2^{<\P}$ such that,
for every HTP-generic $V\in\U_\sigma$, $R_V$ admits a diophantine model of $\Z$
via these three polynomials.
\end{enumerate}
If (2) holds, then $\P\in\D$ (i.e., $\Q$ admits a diophantine model of $\Z$).
\end{thm}
\begin{pf}
For each triple $\hvec=(h,h_+,h_{\times})$ of polynomials of appropriate lengths over $\Z$,
we set $\D_{\hvec}$ to be the class of those $W$ for which $\hvec$
defines a diophantine model of $\Z$ within $R_W$.  If every one of these
classes is nowhere dense, then their countable union $\D$ is meager.

Now suppose that (1) is false, so some class $\D_{\hvec}$ fails to be nowhere dense.
Then there must be
a string $\sigma$ such that $\U_\sigma\subseteq \cl{\D_{\hvec}}$.  Using
this $\sigma$ and this $\hvec$, we now show that in fact all of $\U_\sigma$
is contained within $\D_{\hvec}$.  Let $R_0=R_{\sigma^{-1}(1)}$ and
$R_1=R_{\P-\sigma^{-1}(0)}$ be the largest and smallest subrings
(under $\subseteq$) in $\U_\sigma$, so $R_0$ is finitely generated
and $R_1$ is semilocal.

Fix a single $W\supset\sigma$ with $W\in\D_{\hvec}$, and fix the tuples
$\xvec_0$ and $\xvec_1$ from $R_W$ which represent the elements
$0$ and $1$ in the diophantine model of $\Z$ defined by $\hvec$ in $R_W$.
It follows that $h_{\times}(\xvec_0,\xvec_0,\xvec_0,\Yvec)\in\HTP(R_W)$
and $h_{\times}(\xvec_1,\xvec_1,\xvec_1,\Yvec)\in\HTP(R_W)$.  Now
if any other tuple $\xvec$ from $R_1$ had $h(\xvec,\Yvec)\in\HTP(R_1)$
and $h_{\times}(\xvec,\xvec,\xvec,\Yvec)\in\HTP(R_1)$, then we could set
$\tau=\sigma\hat{~}111\cdots 1$ to contain enough primes that $R_{\tau^{-1}(1)}$
would contain $\xvec$, $\xvec_0$, and $\xvec_1$.  This would mean that
$\hvec$ could not define a diophantine model of $\Z$ in any $R_V$ with
$V\in\U_\tau$, contrary to hypothesis.  Therefore, no other $\xvec$ from $R_1$
can do this.  Now suppose that $\xvec_0$ does \emph{not} lie within $R_0$.
In this case, some extension $\rho=\sigma\hat{~}000\cdots 0$ would exclude
enough primes to ensure that $\xvec_0$ does not lie in $R_{\P-\rho^{-1}(0)}$,
and then no extension of $\rho$ would admit a diophantine model via $\hvec$,
since no other tuple with the right properties lies in $R_1$.  Again, this would
contradict our hypothesis that $\U_\sigma\subseteq \cl{\D_{\hvec}}$,
since $\D_{\hvec}\cap\U_{\rho}$ would be empty,
and so in fact $\xvec_0$ lies in $R_0$, and similarly so does $\xvec_1$.

Now one proceeds by induction on the subsequent elements of the diophantine
model in $R_1$.  Some tuple $\xvec_2$ from $R_W$ must satisfy
$h(\xvec_2,\Yvec)\in\HTP(R_W)$ and
$h_+(\xvec_1,\xvec_1,\xvec_2,\Yvec)\in\HTP(R_W)$, and by the same arguments
as above, we see that $\xvec_2$ is the only tuple in $R_1$ with this property,
and then that $\xvec_2$ actually lies in $R_0$.  Likewise, $\xvec_{-1}$
must satisfy $h(\xvec_{-1},\Yvec)\in\HTP(R_W)$ and
$h_+(\xvec_1,\xvec_{-1},\xvec_0,\Yvec)\in\HTP(R_W)$, and again this
forces $\xvec_{-1}$ to lie in $R_0$ and to be the unique tuple with these
properties in $R_1$.

Continuing this induction, we see that every tuple in the domain of the diophantine
modle of $\Z$ in $R_W$ actually lies in $R_0$, and hence in every $R_W$
with $W\in\U_\sigma$; and moreover that these are the only tuples $\xvec$
in $R_1$ for which $h(\xvec,\Yvec)\in\HTP(R_1)$.  Likewise, if some
$\xvec_m$, $\xvec_n$ and $\xvec_p$ (representing $m$, $n$, and $p$
in the diophantine model) satisfy $h_+(\xvec_m,\xvec_n,\xvec_p,\Yvec)\in\HTP(R_1)$,
then for some $k$, $\tau=\sigma\hat{~}1^k$ is long enough to ensure
that every $W$ extending $\tau$ must have
$h_+(\xvec_m,\xvec_n,\xvec_p,\Yvec)\in\HTP(R_W)$.  But some such $W$
lies in $\D_{\hvec}$, so we must have $m+n=p$.  The same works for $h_{\times}$,
and so $\hvec$ defines a diophantine model of $\Z$ in $R_1$ specifically.

Now it is not clear whether $\hvec$ defines a diophantine model in the
subring $R_0$ (which, being finitely generated, lies in $\B$).
The domain elements of the model in $R_1$
all lie in $R_0$, but the witnesses might not.  However, suppose that
$V\in\U_\sigma$ is HTP-generic, and fix any domain element $\xvec$.
Let $\tau=V\res m$, for any $m\geq |\sigma|$.  Then some $U\supseteq\tau$
lies in $\D_{\hvec}$, and so some extension of $\tau$ yields a solution
to $h(\xvec,\Yvec)$.  Since $V$ is HTP-generic (that is, $V\notin\B$),
this forces $h(\xvec,\Yvec)\in\HTP(R_V)$.
Likewise, for each fact coded by $h_+$ or $h_{\times}$ about domain
elements of the model, some extension of $V\res m$ must yield
a witness to that fact, and therefore $R_V$ itself contains such a witness.  Therefore, $\hvec$
also defines this same diophantine model in \emph{every} HTP-generic
subring $R_V$ with $V\in\U_\sigma$, as required by (2).

Cases (1) and (2) of the theorem cannot both hold, because under (2),
$\U_\sigma\cap\Bbar$ would be a nonmeager subset of $\D$.  Moreover,
the $1$-reduction
$\HTP(R_1)\leq_1\HTPQ$ given in \cite[Proposition 5.4]{EMPS15}
has sufficient uniformity that the images of $h$, $h_+$, and $h_{\times}$
under this reduction define a diophantine model of $\Z$ inside $\Q$.
(Specifically, $h(\Xvec,\Yvec)$ maps to the sum of $h^2$ with several other squares
of polynomials in such a way as to guarantee that all solutions
use values from $R_1$ for the variables $\Xvec$ and $\Yvec$;
likewise with $h_+$ and $h_{\times}$.)
This proves the final statement of the theorem.
\qed\end{pf}

\comment{
NEXT QUESTION:  we should be able to go from Case (2) to a (different) triple
$\hvec$ which works the same with the empty string in place of $\sigma$,
i.e., which defines a diophantine model in every last subring $R_W$ of $\Q$.
Proposition \ref{prop:semilocal} and Julia Robinson's theorem about finitely generated
subrings of $\Q$ seem like the keys.  How to do this?  Is it obvious?
}

Now we continue with the question of existential definability of the integers.

\begin{defn}
\label{defn:existential}
In a ring $R$, a polynomial $g\in\Z[X,\Yvec]$
\emph{existentially defines $\Z$} if, for every $q\in R$,
$$ q\in\Z\iff g(q,\Yvec)\in\HTP(R).$$
$\Z$ is \emph{existentially definable in $R$} if such
a polynomial $g$ exists.
\end{defn}

A ring in which $\Z$ is existentially definable must admit
a very simple diophantine model of $\Z$, given by the polynomial $g$
along with $h_+=X_1+X_2-X_3$ and $h_{\times}=X_1X_2-X_3$.
The question of definability of $\Z$ in the field $\Q$ was originally
answered by Julia Robinson (see \cite{R49}), who gave a $\Pi_4$
definition.  Subsequent work by Poonen \cite{P09}
and then Koenigsmann \cite{K16} has resulted in a $\Pi_1$
definition of $\Z$ in $\Q$, but it remains unknown
whether any existential formula defines $\Z$ there.

\begin{defn}
\label{defn:Eclass}
$\E$ is the class of subrings of $\Q$ in which $\Z$ is existentially definable:
$$ \E = \set{W\subseteq\P}{\Z\text{~is existentially definable in~}R_W}.$$
\end{defn}

We now address the question of the size of the class $\E$.  As with $\D$,
we show $\E$ to be either very large or very small, in the sense of Baire category.

\begin{thm}
\label{thm:existential}
The following are equivalent.
\begin{enumerate}
\item
The class $\E$ is not meager.
\item
There is a $\sigma\in 2^{<\P}$,
and a single polynomial $g$ which existentially defines $\Z$ in
all HTP-generic subrings $R_V$ with $V\in\U_{\sigma}$.
\item
$\P\in\E$; that is, $\Z$ is existentially definable in $\Q$.
\item
There is a single existential formula which defines $\Z$
in every subring of $\Q$.
\end{enumerate}
\end{thm}
\begin{pf}
The proof that $(1)\implies(2)\implies(3)$
proceeds along the same lines as that of Theorem \ref{thm:diophantine},
with $\E_g$ as the class of those $W$ for which the polynomial $g$
existentially defines $\Z$ within $R_W$.  If every one of these
classes is nowhere dense, then their countable union $\E$ is meager.
Otherwise one proves (2), and from that (3), by a simplification of the
same method as before, wth no induction required.
To see that (3) implies (4), notice that if $\Z$ is defined in $\Q$
by the formula $\exists\Yvec~f(X,\Yvec)=0$, and $d$ is the
total degree of $f$, then the formula
$$\exists\Yvec\exists Z~[Z^d\cdot f\left(X,\frac{Y_1}Z,\ldots,\frac{Y_n}Z\right)=0~\&~Z>0]$$
defines $\Z$ in $R_W$; this is the same trick we used in Section \ref{sec:intro}
to reduce $\HTPQ$ to $\HTP(R_W)$.
\qed\end{pf}

It is possible to turn Theorem \ref{thm:diophantine} into an equivalence
analogous to that in Theorem \ref{thm:existential}, with the third condition
stating that $\P\in\D$.  As far as we know, however, it is necessary
to consider diophantine \emph{interpretations} in subrings $R_W$,
rather than diophantine models, in order to accomplish this.  The distinction
is simply that in a diophantine interpretation of $\Z$ in $R_W$, the domain is
allowed to be a diophantine subset of $\R_W^n$ modulo
a diophantine equivalence relation, with operations (still
defined by diophantine formulas) that respect this equivalence relation.
It is readily seen that a diophantine interpretation of the ring $\Z$
in $\Q$ yields a uniform diophantine interpretation of $\Z$ in
every $R_W$, by the same device as in the proof of Theorem \ref{thm:existential}:
each domain element $x$ of the interpretation in $\Q$ may be represented
in any $R_W$ by those pairs $(y,z)$ with $xz=y$ and $z>0$, modulo
the obvious equivalence relation.

\section{Measure Theory}
\label{sec:measuretheory}

Normally there is a strong connection between measure theory and
Baire category theory.  Each defines a certain $\Sigma$-ideal of sets
to be ``small'':  the sets of measure $0$, and the meager sets, respectively.
In Cantor space, as in their original model $\mathbb R$,
neither of these two properties is strong enough to imply the other,
but empirically they appear closely connected:  sets of measure $0$
are very often meager, and vice versa, unless the sets are specifically selected
to avoid this.  (One difference was mentioned above, in the context of
Lemma \ref{lemma:uniform}.)

Our results in this article rely heavily on the simple Lemma \ref{lemma:boundarymeager},
which stated that the boundary set $\B(f)$ of a polynomial $f$ must be nowhere dense.
Most of our subsequent results have measure-theoretic analogues which
would go through fairly easily, provided that boundary sets $\B(f)$
also have measure $0$.  However, determining the measure of
the boundary set of a polynomial appears to be a nontrivial problem.
It is unknown whether there exists any polynomial $f$ for which
the measure $\mu(\B(f)) >0$.  Moreover, if such an $f$ exists,
it is unclear what other constraints on the real number $\mu(\B(f))$
exist, apart from the computability-theoretic upper bound
given by its definition as $\mu(\B(f))$.  Could such a number be transcendental?
Could it be a noncomputable real number?  If not, is there an algorithm
for computing $\mu(\B(f))$ uniformly in $f$?  All of these
appear to be challenging questions, often with a more number-theoretic
flavor than most of this article.  If they can be resolved,
then it may be possible to determine whether or not Hilbert's Tenth Problem
on subrings of $\Q$ has measure-theoretic zero-one laws similar
to those proven here for Baire category.

\section{Acknowledgments}

The author was partially supported by Grant \# DMS -- 1362206
from the National Science Foundation, and by several grants
from the PSC-CUNY Research Award Program.
This work grew out of research initiated at a workshop held at the
American Institute of Mathematics and continued at a workshop held at the
Institute for Mathematical Sciences of the National University of Singapore.
Conversations with Bjorn Poonen and Alexandra Shlapentokh have been
very helpful in the creation of this article.

\parbox{4.7in}{
{\sc
\noindent
Department of Mathematics \hfill \\
\hspace*{.1in}  Queens College -- C.U.N.Y. \hfill \\
\hspace*{.2in}  65-30 Kissena Blvd. \hfill \\
\hspace*{.3in}  Queens, New York  11367 U.S.A. \hfill \\
Ph.D. Programs in Mathematics \& Computer Science \hfill \\
\hspace*{.1in}  C.U.N.Y.\ Graduate Center\hfill \\
\hspace*{.2in}  365 Fifth Avenue \hfill \\
\hspace*{.3in}  New York, New York  10016 U.S.A. \hfill}\\
\medskip
\hspace*{.045in} {\it E-mail: }
\texttt{Russell.Miller\at {qc.cuny.edu} }\hfill \\
}

\end{document}